\documentclass{amsart}

\usepackage{amsmath,amsthm,amssymb}
\usepackage{hyperref}

\usepackage{booktabs}
\usepackage{longtable}
\usepackage{array}

\usepackage{float}
\usepackage{graphicx}

\newtheorem{proposition}{Proposition}
\newtheorem{theorem}{Theorem}
\newtheorem{corollary}{Corollary}

\newtheorem{example}{Example}

\title[Binomial AS-Regular Algebras]
{Binomial Artin--Schelter Regular Algebras and Twisted Monoid Algebras}

\author{Vishal Bhatoy}
\address{Department of Mathematics, Carleton University, Ottawa, Ontario, Canada}
\email{vishalbhatoy@cmail.carleton.ca}

\begin{document}

\begin{abstract}
We study four-dimensional binomial
Artin--Schelter regular algebras and
their connection with twisted monoid algebras.
Starting with $23$ examples up to permutation
of generators, we investigate the associated
monoids of I-type and the structure of their
groups of fractions.
We determine the corresponding finite
permutation groups and compute the second
group cohomology $H^2(G(M),k^\times)$.
These computations allow us to construct
explicit families of four-dimensional
Koszul Artin--Schelter regular algebras
through cocycle twists. We illustrate this construction for
some of these algebras.
\end{abstract}

\maketitle

\section{Introduction}

In this paper, we study binomial
Artin--Schelter (AS) regular algebras
through their connections with twisted
monoid algebras. Our main interest is
in the monoids arising from binomial
relations, their groups of fractions
and the cocycles that determine the
coefficients in these relations.

Let $A$ be a connected graded domain
equipped with a faithful grading by
a group $G$ such that
\[
A=\bigoplus_{g\in G}A_g,
\qquad
\dim_k A_g\leq 1.
\]
The support of this grading forms a
submonoid $M$ of $G$, and multiplication
in $A$ determines a $2$-cocycle on $M$.
Consequently, $A$ can be described as
a twisted monoid algebra
\[
A\cong k_\sigma[M].
\]
The monoid $M$ records the relations
among the degrees of the generators,
while the cocycle $\sigma$ determines
the coefficients in their products.

This motivates the study of
coefficient-free relations obtained
by replacing the nonzero coefficients
in the binomial relations by $1$.
When these relations present $M$,
the resulting algebra is the monoid
algebra $k[M]$. We investigate when
the original algebra can be recovered
as a cocycle twist of $k[M]$.

If the monoid $2$-cocycle extends
to a group $2$-cocycle and the
group grading is compatible with
the connected $\mathbb N$-grading,
then cocycle twisting preserves
Koszul AS-regularity and global
dimension
\cite[Theorem~4.2.1]{HuangEtAl2024}.

An important class of examples comes
from monoids of I-type, introduced
by Gateva-Ivanova and Van den Bergh
\cite{GatevaIvanovaVanDenBergh1998}.
These monoids are closely related
to finite involutive nondegenerate
set-theoretic solutions of the
Yang--Baxter equation
\cite{EtingofSchedlerSoloviev1999}.
Jespers and Okni\'nski showed that
the groups of fractions of monoids
of I-type on $n$ generators admit
embeddings into
$\mathbb Z^n\rtimes S_n$
and contain finite-index free abelian
subgroups of rank $n$
\cite{JespersOkninski2005}.
In particular, these groups are
solvable Bieberbach groups of
Hirsch length $n$.

We focus particularly on
four-dimensional Koszul AS-regular
algebras defined by four generators
and six quadratic binomial relations.
We started with a list of $168$ binomial
algebras, which is available on
\href{https://github.com/ColinIngalls/invariants-of-regular-algebras/blob/master/KoszulCohomologyAndKodairaSpencerMap%20-%20GAPCode/CodeFiles/binom4.m2.txt}
{GitHub}.
For each algebra, we computed its
superpotential in the sense of Mori
and Smith \cite{Mori_Izuru_Smith}.
We then considered the natural action
of the symmetric group $S_4$ on the
variables $x_0,x_1,x_2,x_3$.
Two superpotentials were considered
equivalent if one could be obtained
from the other by permuting the
variables. This gave $23$
representative superpotentials
under the action of $S_4$.

For these $23$ representatives,
we study the associated coefficient-free
monoids, their groups of fractions,
and the finite permutation groups
arising from the corresponding
Yang--Baxter solutions.
We also compare these groups with
known examples in the literature,
including those studied by Jespers
and Okni\'nski
\cite{JespersOkninski2005}
and Buzaglo and Rogalski
\cite{BuzagloRogalski2026}.

We also compute the second group
cohomology $H^2(G(M),k^\times)$
for the associated groups of fractions.
This group classifies group
$2$-cocycles up to coboundaries
and provides information about
cocycle twists of the corresponding
monoid algebras.

The paper is organized as follows.
In Section ~\ref{prelim}, we recall the necessary
background on group gradings, cocycles,
binomial AS-regular algebras, monoids
of I-type and their groups of fractions.
In Section~\ref{binomasregulaandtwistedmonoid}, we explain the connection
between binomial AS-regular algebras
and twisted monoid algebras, including
the role of cocycles in determining
the coefficients of binomial relations.
In Section~\ref{groupoffractions}, we study monoids of I-type
and the structure of their groups
of fractions.
In Section~\ref{sec:low-dimensional-examples}, we illustrate these
constructions through examples in
dimensions $2$ and $3$.
In Section~\ref{sec:four-dimensional-examples}, we investigate the
$23$ four-dimensional examples and
determine their associated permutation
groups.
In Section~\ref{sec:computationscocycles}, we compute the second
cohomology groups $H^2(G(M),k^\times)$
for the examples considered. In Section~\ref{sec:explicit-cocycle-twists},
we use these computations to construct
explicit cocycle twists of some algebras.
Finally, the appendix~\ref{app:binomial-algebras} contains the
defining relations of the
$23$ four-dimensional binomial algebras.

\section{Preliminaries}\label{prelim}

Throughout this paper, let $k$ be an
algebraically closed field of
characteristic zero. We recall some definitions
and results that we will use throughout
the paper. We begin with group gradings
and cocycles, followed by binomial
AS-regular algebras, monoids of I-type
and their groups of fractions.

\subsection{Group gradings and cocycles}

Let $G$ be a group.
An algebra $A$ is called $G$-graded if
\[
A=\bigoplus_{g\in G}A_g,
\qquad
A_gA_h\subseteq A_{gh}
\]
for all $g,h\in G$.

The grading is called faithful if
the support
\[
\operatorname{supp}(A)
=
\{g\in G\mid A_g\neq 0\}
\]
generates $G$ as a group.
If $A$ is a domain, then
$\operatorname{supp}(A)$ is a
submonoid of $G$.

A normalized $2$-cocycle on $G$
is a function
\[
\sigma:G\times G\longrightarrow k^\times
\]
satisfying
$\sigma(g,h)\sigma(gh,\ell)
=
\sigma(h,\ell)\sigma(g,h\ell)$
for all $g,h,\ell\in G$, and
\[
\sigma(1,g)=\sigma(g,1)=1.
\]

Two normalized $2$-cocycles
$\sigma$ and $\tau$ are called
\textit{cohomologous} if there exists a
function $f:G\to k^\times$,
with $f(1)=1$, such that
\[
\tau(g,h)
=
\sigma(g,h)
\frac{f(g)f(h)}{f(gh)}
\]
for all $g,h\in G$.
The equivalence classes of
normalized $2$-cocycles form
the second cohomology group
$H^2(G,k^\times)$, where $G$
acts trivially on $k^\times$.

A normalized $2$-cocycle on a monoid
$M$ is defined in the same way,
with $g,h,\ell\in M$.

Given a $G$-graded algebra $A$
and a normalized $2$-cocycle
$\sigma$ on $G$, we define a new
multiplication on homogeneous
elements by
\[
a*_\sigma b
=
\sigma(g,h)ab,
\qquad
a\in A_g,\quad b\in A_h.
\]
The cocycle identity ensures that
this multiplication is associative.
The resulting algebra, denoted by
$A^\sigma$, is called a
cocycle twist of $A$.
It has the same underlying graded
vector space as $A$.

Cocycle twists are closely related
to Zhang twists
\cite{Davies2017}.
Under suitable assumptions,
cocycle twisting preserves
Koszul AS-regularity
\cite{HuangEtAl2024}.

\subsection{Binomial Artin--Schelter regular algebras}

Recall that a connected
$\mathbb N$-graded algebra has the form
\[
A=\bigoplus_{n\geq 0}A_n,
\qquad A_0=k.
\]

A connected graded algebra is called
Artin--Schelter regular
\cite{Artin1987} if it has
finite global dimension, finite
Gelfand--Kirillov dimension, and
satisfies the AS-Gorenstein condition.
More precisely, if its global dimension
is $d$, then
\[
\operatorname{Ext}_A^i(k,A)=0
\qquad (i\neq d),
\]
and
\[
\operatorname{Ext}_A^d(k,A)\cong k(\ell)
\]
for some integer $\ell$, where
the Ext groups are taken in the
graded category.

An algebra is called Koszul
\cite{polishchuk2005quadratic}
if the trivial module $k$ has
a graded free resolution in which
the generators at homological
degree $i$ occur in internal degree $i$.

In this paper, we consider Koszul
Artin--Schelter regular algebras
generated in degree one and presented
by quadratic binomial relations
of the form
\[
x_ax_b=q\,x_cx_e,
\qquad q\in k^\times.
\]

We are particularly interested
in the coefficient-free relations
obtained by setting the nonzero
coefficients equal to $1$.
These relations define a monoid
\[
M=\langle x_1,\ldots,x_n
\mid x_ax_b=x_cx_e\rangle^+,
\]
whose group of fractions we study
when it exists.

\subsection{Monoids of I-type and Yang--Baxter solutions}

Monoids of I-type were introduced by
Gateva-Ivanova and Van den Bergh
\cite{GatevaIvanovaVanDenBergh1998}.
They are closely related to
set-theoretic solutions of the
Yang--Baxter equation
\cite{EtingofSchedlerSoloviev1999}.

Let $X$ be a finite set.
A set-theoretic solution of the
Yang--Baxter equation is a map
\[
r:X\times X\longrightarrow X\times X
\]
satisfying
\[
(r\times\operatorname{id})
(\operatorname{id}\times r)
(r\times\operatorname{id})
=
(\operatorname{id}\times r)
(r\times\operatorname{id})
(\operatorname{id}\times r).
\]
We follow the terminology of
\cite[Section~1]{JespersOkninski2005}.

For each $i,j\in X$, write
\[
r(i,j)=\bigl(L_i(j),R_j(i)\bigr),
\]
where
\[
L_i(j)=\pi_1(r(i,j)),
\qquad
R_j(i)=\pi_2(r(i,j)),
\]
and $\pi_1,\pi_2:X\times X\to X$
are the natural projections.
The maps $L_i:X\to X$ and
$R_j:X\to X$ are called the
left and right actions, respectively.

The solution is called involutive if
$r^2=\operatorname{id}$, and
nondegenerate if every $L_i$ and
$R_j$ is a bijection of $X$.

The left actions generate a finite
permutation group
\[
F=\langle L_i\mid i\in X\rangle
\leq\operatorname{Sym}(X).
\]

The structure monoid associated with
$(X,r)$ is
\[
M(X,r)=
\left\langle
x_i\ (i\in X)
\ \middle|\
x_ix_j=x_{L_i(j)}x_{R_j(i)}
\right\rangle^+.
\]
The relationship between these
monoids and Yang--Baxter solutions
is described in
\cite[Theorem~1.1]{JespersOkninski2005}.

\begin{theorem}[\cite{JespersOkninski2005}, Theorem~2.2 and
Corollary~2.4]
\label{thm:itype-structure}
Let $(X,r)$ be a finite involutive
nondegenerate set-theoretic solution
of the Yang--Baxter equation, with
$|X|=d$, and let $M=M(X,r)$.
Then $M$ is a monoid of I-type
and has a group of fractions $G(M)$.
Moreover, there is an embedding
\[
G(M)\hookrightarrow
\mathbb Z^d\rtimes S_d
\]
such that the restriction of the
set-theoretic projection
\[
\mathbb Z^d\rtimes S_d
\longrightarrow\mathbb Z^d
\]
to the image of $G(M)$ is bijective.

The image of $G(M)$ under the
projection onto $S_d$ is the
permutation group
\[
F=\langle L_i\mid i\in X\rangle.
\]
Furthermore, $G(M)$ is a Bieberbach
group and fits into an exact sequence
\[
1\longrightarrow N
\longrightarrow G(M)
\longrightarrow F
\longrightarrow 1,
\]
where $N$ is a free abelian group
of rank $d$. In particular,
$G(M)$ is solvable,
virtually abelian, and has Hirsch length
\[
h(G(M))=d.
\]
\end{theorem}

\subsection{Groups of fractions and Hirsch length}

A group is called \textit{virtually abelian}
if it contains an abelian subgroup of
finite index. Similarly, it is called
\textit{virtually nilpotent} if it contains
a nilpotent subgroup of finite index.

A group $G$ is called \textit{polycyclic}
if it has a subnormal series
\[
1=G_0\triangleleft G_1
\triangleleft\cdots\triangleleft G_m=G
\]
such that each quotient $G_i/G_{i-1}$
is cyclic.

The \textit{Hirsch length} $h(G)$
is the number of infinite cyclic
factors in such a series.
It is independent of the chosen
series
\cite[Theorem~5.4.13]{Robinson1996}.

We recall the following standard
property of Hirsch length. 

\begin{proposition}[\cite{Robinson1996}, Section~5.4,
Exercise~9]
\label{prop:hirsch-additivity}
Let
\[
1\longrightarrow N
\longrightarrow G
\longrightarrow Q
\longrightarrow 1
\]
be an exact sequence of polycyclic
groups. Then
\[
h(G)=h(N)+h(Q).
\]
In particular, if $Q$ is finite,
then
\[
h(G)=h(N).
\]
\end{proposition}

For a monoid of I-type on $d$
generators, Theorem~\ref{thm:itype-structure}
gives an exact sequence
\[
1\longrightarrow N
\longrightarrow G(M)
\longrightarrow F
\longrightarrow1,
\]
where $N\cong\mathbb Z^d$ and $F$
is finite.
Applying Proposition~\ref{prop:hirsch-additivity},
we obtain
\[
h(G(M))=h(N)=d.
\]

These properties will be useful
when studying the groups arising
from coefficient-free binomial
relations in dimensions $2$, $3$
and $4$.

\section{Binomial Artin--Schelter regular algebras as twisted monoid algebras}\label{binomasregulaandtwistedmonoid}

We explain how the monoid associated with
a binomial algebra records the relations
among the degrees of its generators, while
a cocycle determines the coefficients
appearing in these relations.

\subsection{Twisted monoid algebras}

Let $A$ be a connected graded domain
with a faithful grading
\[
A=\bigoplus_{g\in G}A_g
\]
by a group $G$, and suppose that
$\dim_k A_g\leq 1$ for every $g\in G$.

Let
\[
M=\{g\in G\mid A_g\neq 0\}
\]
be the support of the grading.
Since $A$ is a domain, $M$ is a
submonoid of $G$, and faithfulness
implies that $M$ generates $G$.

\begin{proposition}
\label{prop:twisted-monoid}
There is a normalized $2$-cocycle
\[
\sigma:M\times M\longrightarrow k^\times
\]
such that
\[
A\cong k_\sigma[M]
\]
as $M$-graded algebras.
\end{proposition}

\begin{proof}
For each $m\in M$, choose a nonzero
element $u_m\in A_m$, with $u_1=1$.
Since $A_m$ is one-dimensional,
we have $A_m=ku_m$.

For $m,n\in M$, the product $u_mu_n$
is nonzero and belongs to $A_{mn}$.
Thus there is a unique scalar
$\sigma(m,n)\in k^\times$ such that
\[
u_mu_n=\sigma(m,n)u_{mn}.
\]

Associativity gives
\[
\sigma(m,n)\sigma(mn,p)
=
\sigma(n,p)\sigma(m,np),
\]
while $u_1=1$ implies
\[
\sigma(1,m)=\sigma(m,1)=1.
\]
Hence $\sigma$ is a normalized
$2$-cocycle, and the correspondence
$e_m\mapsto u_m$ defines an
$M$-graded algebra isomorphism
\[
k_\sigma[M]\cong A.
\]
\end{proof}

The cocycle $\sigma$ is defined on
$M$. When it extends to a group
$2$-cocycle on $G$, it also determines
a class in $H^2(G,k^\times)$.

\subsection{Binomial relations and their coefficients}

Suppose that $A$ is generated by
nonzero $G$-homogeneous elements
$x_1,\ldots,x_d$ and has quadratic
binomial relations
\[
m_i-q_i n_i,
\qquad q_i\in k^\times.
\]
Write $g_j=\deg_G(x_j)$.

For a relation
\[
x_ax_b=q_i x_cx_e,
\]
homogeneity gives
\[
g_ag_b=g_cg_e.
\]

Using the basis from
Proposition~\ref{prop:twisted-monoid},
write
\[
x_j=c_j u_{g_j},
\qquad c_j\in k^\times.
\]
Then
\[
x_ax_b
=
c_ac_b\,\sigma(g_a,g_b)u_{g_ag_b},
\]
and
\[
x_cx_e
=
c_cc_e\,\sigma(g_c,g_e)u_{g_cg_e}.
\]
Consequently,
\[
q_i=
\frac{c_ac_b\,\sigma(g_a,g_b)}
{c_cc_e\,\sigma(g_c,g_e)}.
\]

Thus the relations among the degrees
are recorded by $M$, while the
coefficients arise from the cocycle
$\sigma$ together with the choice
of homogeneous generators.

\subsection{Removing the coefficients}

Suppose
\[
A=
k\langle x_1,\ldots,x_d\rangle/
(m_i-q_i n_i)_{i=1}^r,
\qquad q_i\in k^\times,
\]
and assume that the coefficient-free
relations present the support monoid
\[
M=
\left\langle
g_1,\ldots,g_d
\ \middle|\
m_i=n_i,\quad 1\leq i\leq r
\right\rangle^+.
\]

Define
\[
A_{\mathrm{cf}}
=
k\langle x_1,\ldots,x_d\rangle/
(m_i-n_i)_{i=1}^r.
\]
By the presentation of $M$,
\[
A_{\mathrm{cf}}\cong k[M].
\]
On the other hand,
Proposition~\ref{prop:twisted-monoid}
gives
\[
A\cong k_\sigma[M].
\]
Thus $A$ is obtained by twisting
the multiplication of $k[M]$,
and
\[
A\cong(k[M])^\sigma,
\qquad
k[M]\cong A^{\sigma^{-1}},
\]
where the twists are taken with
respect to the $M$-grading.

We next consider when this
construction preserves
Koszul AS-regularity.

\begin{proposition}
\label{prop:coefficient-removal}
Suppose that $A$ is a Koszul
AS-regular domain and that the
monoid cocycle $\sigma$ extends
to a normalized group $2$-cocycle
on $G$. Assume that the
$G$-grading is compatible with
the connected $\mathbb N$-grading,
in the sense that each $A_g$
is a graded subspace.

Then $A^{\sigma^{-1}}$ is Koszul
AS-regular of the same global
dimension as $A$. Consequently,
\[
A_{\mathrm{cf}}\cong k[M]
\cong A^{\sigma^{-1}}
\]
is Koszul AS-regular of the same
global dimension as $A$.
\end{proposition}

\begin{proof}
The $G$-grading defines a right
coaction of the group algebra $kG$
on $A$. By the compatibility
assumption, this coaction preserves
the connected $\mathbb N$-grading.

The inverse group cocycle
$\sigma^{-1}$ defines a Hopf
$2$-cocycle on $kG$.
By
\cite[Theorem~4.2.1]{HuangEtAl2024},
the corresponding cocycle twist
preserves Koszul AS-regularity
and global dimension.

Therefore $A^{\sigma^{-1}}$ is
Koszul AS-regular of the same
global dimension as $A$.
Since
\[
A^{\sigma^{-1}}\cong k[M]
\cong A_{\mathrm{cf}},
\]
the result follows.
\end{proof}

\section{Monoids and their groups of fractions}\label{groupoffractions}

We study the monoids arising from
coefficient-free quadratic binomial
relations and their associated groups.

Our main examples come from
four-dimensional Koszul AS-regular
algebras with four generators and
six quadratic binomial relations.
We consider monoids of the form
\[
M=
\left\langle
g_1,g_2,g_3,g_4
\ \middle|\
m_i=n_i,\quad 1\leq i\leq 6
\right\rangle^+,
\]
where $m_i$ and $n_i$ are words
of length two.

The group defined by the same
presentation is
\[
U(M)=
\left\langle
g_1,g_2,g_3,g_4
\ \middle|\
m_i=n_i,\quad 1\leq i\leq 6
\right\rangle.
\]
This is the universal group of $M$.
When $M$ is cancellative and satisfies
the Ore condition, it embeds in its
group of fractions, which we denote
by $G(M)$. In this case,
\[
G(M)\cong U(M).
\]

\subsection{Binomial monoids of I-type}

We first recall a criterion for
recognizing binomial monoids of I-type.

\begin{theorem}[\cite{GatevaIvanovaVanDenBergh1998},
Theorems~1.1 and~1.3]
\label{thm:binomial-monoid-itype}
Let $M$ be a monoid generated by
$g_1,\ldots,g_d$ with
$\binom{d}{2}$ quadratic relations
\[
g_jg_i=g_{i'}g_{j'},
\qquad j>i.
\]
Suppose that, after a suitable
ordering of the generators, the
relations satisfy the following
conditions:
\begin{enumerate}
\item Each relation has
$i'<j'$ and $i'<j$.

\item As the relations vary,
every pair $(i',j')$ with
$i'<j'$ occurs exactly once.

\item The overlaps $g_kg_jg_i$,
where $k>j>i$, give rise to
no new relations.
\end{enumerate}
Then $M$ is a monoid of I-type.
\end{theorem}

The preceding theorem concerns
binomial monoids with square-free
quadratic relations. More general
monoids of I-type may also have
relations involving squares.
For these monoids, we use the
characterization in
\cite[Theorem~1.1]{JespersOkninski2005}.

\subsection{Structure of the groups of fractions}

The following result records some
properties of monoids of I-type
and their groups of fractions.

\begin{proposition}[\cite{JespersOkninski2005}, Corollaries~8.2.5
and~8.2.7]
\label{prop:itype-fractions}
Let $M$ be a monoid of I-type
on $d$ generators. Then:
\begin{enumerate}
\item $M$ is cancellative and
admits a group of fractions $G(M)$.

\item $G(M)$ is torsion-free,
solvable and virtually abelian.

\item There is an exact sequence
\[
1\longrightarrow\mathbb Z^d
\longrightarrow G(M)
\longrightarrow F
\longrightarrow1,
\]
where $F$ is a finite
permutation group.

\item $G(M)$ is a Bieberbach
group with Hirsch length $d$.
\end{enumerate}
\end{proposition}

In particular, the group of fractions
of a monoid of I-type is virtually
nilpotent and polycyclic.

\begin{corollary}
\label{cor:non-itype}
Let $M$ be a finitely generated
monoid admitting a group of
fractions $G(M)$. If $G(M)$ is
not virtually abelian, then
$M$ is not a monoid of I-type.
\end{corollary}

\begin{proof}
By
Proposition~\ref{prop:itype-fractions},
the group of fractions of a
monoid of I-type is virtually
abelian. The assertion follows
by contraposition.
\end{proof}

\subsection{Algebras of monoids of I-type}

The structure of a monoid of I-type
also imposes strong restrictions
on its monoid algebra.

\begin{proposition}[\cite{JespersOkninski2007},
Theorem~8.5.1]
\label{prop:itype-algebra}
Let $M$ be a monoid of I-type
on $d$ generators. Then $k[M]$
is a Noetherian domain satisfying
a polynomial identity. Moreover,
\[
\operatorname{GKdim} k[M]=d.
\]
\end{proposition}

Furthermore, Gateva-Ivanova and
Van den Bergh
\cite[Theorem~1.4]{GatevaIvanovaVanDenBergh1998}
showed that if $M$ is a monoid
of I-type and $\sigma$ is a
monoid $2$-cocycle, then the
twisted monoid algebra
$k_\sigma[M]$ is Koszul,
Noetherian and has finite
global dimension.

For our four-dimensional examples,
we investigate whether the
coefficient-free monoids satisfy
the I-type conditions. We also
study the structure of their
groups of fractions, including
solvability, virtual nilpotence
and Hirsch length.

In the examples where a group
of fractions is identified, we
subsequently study
\[
H^2(G(M),k^\times)
\]
to understand the possible
group cocycles associated with
the binomial algebras.

\section{Low-dimensional examples}\label{sec:low-dimensional-examples}

We begin with examples in dimensions
$2$ and $3$, illustrating how the
coefficient-free relations determine
the monoids and their groups of fractions.

In dimension $3$, we compare our
calculations with the classification
of nonabelian group gradings of
three-dimensional Artin--Schelter regular algebras
obtained by Buzaglo and Rogalski
\cite[Theorem~3.4 and Table~1]{BuzagloRogalski2026}.

\subsection{Dimension \texorpdfstring{$2$}{2}}

We consider the two monoids of I-type
on two generators.

\begin{example}
Consider the polynomial ring
\[
A_{\mathrm{cf}}
=
k\langle x,y\rangle/(xy-yx).
\]
Its monoid is
\[
M=\langle x,y\mid xy=yx\rangle^+
\cong\mathbb N^2.
\]
The associated Yang--Baxter solution
is the flip map, with trivial left
actions. Thus
\[
F=1,\qquad
G(M)\cong\mathbb Z^2,
\qquad h(G(M))=2.
\]
\end{example}

\begin{example}
Consider the coefficient-free algebra
\[
A_{\mathrm{cf}}
=
k\langle x,y\rangle/(x^2-y^2),
\]
with associated monoid
\[
M=\langle x,y\mid x^2=y^2\rangle^+.
\]
This monoid is of I-type
\cite[Example~2.5]{JespersOkninski2005}.

The associated left actions are
\[
L_x=L_y=(x\,y),
\]
so
\[
F=\langle(x\,y)\rangle\cong C_2.
\]

The group of fractions has presentation
\[
G(M)=\langle x,y\mid x^2=y^2\rangle.
\]
Setting
\[
a=x,\qquad b=x^{-1}y,
\]
we have $y=ab$. The relation
$x^2=y^2$ becomes
\[
a^2=(ab)^2,
\]
which is equivalent to
\[
aba^{-1}=b^{-1}.
\]
Consequently,
\[
G(M)\cong
\langle a,b\mid aba^{-1}=b^{-1}\rangle,
\]
the Klein bottle group.

The subgroup
\[
N=\langle a^2,b\rangle
\cong\mathbb Z^2
\]
has index $2$. Hence
\[
h(G(M))=2.
\]
\end{example}

Thus the two groups have the same
Hirsch length, although one is
abelian and the other is nonabelian.

\subsection{Dimension \texorpdfstring{$3$}{3}}

We consider a noncommutative example
whose group of fractions is a
semidirect product.

\begin{example}
Consider the coefficient-free algebra
\[
A_{\mathrm{cf}}
=
k\langle x,y,z\rangle/
(xy-yx,\;xz-zy,\;yz-zx).
\]
The associated monoid is
\[
M=
\langle x,y,z\mid
xy=yx,\;xz=zy,\;yz=zx\rangle^+.
\]

The corresponding left actions are
\[
L_x=L_y=1,\qquad L_z=(x\,y).
\]
Therefore
\[
F=\langle(x\,y)\rangle\cong C_2.
\]

The group of fractions is
\[
G(M)=
\langle x,y,z\mid
xy=yx,\;xz=zy,\;yz=zx\rangle.
\]

The last two relations give
\[
z^{-1}xz=y,\qquad
z^{-1}yz=x.
\]
Thus conjugation by $z$ exchanges
$x$ and $y$.

Since $xy=yx$, we obtain
\[
G(M)\cong
\mathbb Z^2\rtimes_\psi\mathbb Z,
\qquad
\psi=
\begin{pmatrix}
0&1\\
1&0
\end{pmatrix}.
\]

This group also appears in the
classification of Buzaglo and Rogalski
\cite[Example~2.18 and Table~1]{BuzagloRogalski2026}.

The subgroup
\[
N=\langle x,y,z^2\rangle
\cong\mathbb Z^3
\]
has index $2$ in $G(M)$.
Consequently,
\[
h(G(M))=3.
\]

\end{example}

These examples illustrate that
the groups of fractions of monoids
of I-type need not be abelian,
even though they contain
finite-index free abelian subgroups.

We next investigate the corresponding
groups for four-dimensional
coefficient-free binomial algebras,
with particular attention to their
structure, Hirsch length and
second cohomology groups.

\section{Four-dimensional examples}
\label{sec:four-dimensional-examples}

We now consider the coefficient-free monoids
associated with the $23$ four-dimensional
binomial Artin--Schelter regular algebras
listed in Appendix~\ref{app:binomial-algebras}.
We study their associated set-theoretic
solutions of the Yang--Baxter equation
and the permutation groups generated
by the left actions.

Let $X=\{0,1,2,3\}$, and write
\[
r(i,j)=\bigl(L_i(j),R_j(i)\bigr)
\]
for the associated solution. The left actions
generate a finite permutation group
\[
F=\langle L_0,L_1,L_2,L_3\rangle
\leq S_4.
\]

For monoids of I-type, the structure theorem
recalled in Proposition~\ref{prop:itype-fractions}
gives an exact sequence
\[
1\longrightarrow\mathbb Z^4
\longrightarrow G(M)
\longrightarrow F
\longrightarrow1.
\]
In particular, $G(M)$ is a torsion-free
Bieberbach group with Hirsch length $4$.

Examples of four-generator monoids of I-type
have also been studied in
\cite[Examples~8.2.14 and~8.4.6]{JespersOkninski2007}.
We consider three algebras from our collection,
illustrating different permutation groups.

\subsection{The commutative algebra}

\begin{example}[Algebra 1]\label{example1}
Consider the coefficient-free relations
\[
\begin{aligned}
x_0x_1&=x_1x_0, &
x_0x_2&=x_2x_0, &
x_0x_3&=x_3x_0,\\
x_1x_2&=x_2x_1, &
x_1x_3&=x_3x_1, &
x_2x_3&=x_3x_2.
\end{aligned}
\]

Let $M$ be the monoid defined by these relations.
Since all generators commute, we have
\[
M\cong\mathbb N^4.
\]

The associated Yang--Baxter solution
is the flip map
\[
r(i,j)=(j,i).
\]
Thus all left actions are trivial:
\[
L_0=L_1=L_2=L_3=1.
\]
Hence the permutation group is
\[
F=\{1\}.
\]

The group of fractions is
\[
G(M)\cong\mathbb Z^4,
\]
which is an abelian Bieberbach group
with Hirsch length $4$.
\end{example}

\subsection{An algebra with cyclic permutation group}

\begin{example}[Algebra 17]\label{example17}
Consider the coefficient-free relations
\[
\begin{aligned}
x_0^2&=x_3^2, &
x_0x_2&=x_1x_3, &
x_0x_3&=x_2x_1,\\
x_1^2&=x_2^2, &
x_1x_2&=x_3x_0, &
x_2x_0&=x_3x_1.
\end{aligned}
\]

Let $M$ be the monoid defined by these relations.
The associated solution has left actions
\[
\begin{aligned}
L_0=L_2&=(0\,3\,2\,1),\\
L_1=L_3&=(0\,1\,2\,3).
\end{aligned}
\]

For example, the relations give
\[
\begin{aligned}
r(0,0)&=(3,3),\\
r(0,1)&=(0,1),\\
r(0,2)&=(1,3),\\
r(0,3)&=(2,1).
\end{aligned}
\]
Therefore
\[
L_0=(0\,3\,2\,1).
\]

The remaining left actions are obtained
similarly. In particular,
\[
L_1=L_3=L_0^{-1},
\qquad L_2=L_0.
\]

Since $L_0$ is a $4$-cycle,
the permutation group is
\[
F=\langle(0\,3\,2\,1)\rangle
\cong C_4.
\]

Since $M$ is a monoid of I-type,
its group of fractions fits into
the exact sequence
\[
1\longrightarrow\mathbb Z^4
\longrightarrow G(M)
\longrightarrow C_4
\longrightarrow1.
\]
In particular, $G(M)$ is a solvable
Bieberbach group with Hirsch length $4$.
\end{example}

\subsection{An algebra with dihedral permutation group}

\begin{example}[Algebra 21]\label{example21}
Consider the coefficient-free relations
\[
\begin{aligned}
x_0^2&=x_3x_1, &
x_0x_1&=x_1x_2, &
x_0x_2&=x_2x_0,\\
x_1^2&=x_3^2, &
x_1x_3&=x_2^2, &
x_2x_3&=x_3x_0.
\end{aligned}
\]

Let $M$ be the monoid defined by these relations.
The associated solution has left actions
\[
\begin{aligned}
L_0&=(0\,3),\\
L_1&=(0\,1\,3\,2),\\
L_2&=(1\,2),\\
L_3&=(0\,2\,3\,1).
\end{aligned}
\]

For example, the relations give
\[
\begin{aligned}
r(1,0)&=(1,0),\\
r(1,1)&=(3,3),\\
r(1,2)&=(0,1),\\
r(1,3)&=(2,2).
\end{aligned}
\]
Therefore
\[
L_1=(0\,1\,3\,2).
\]

The remaining left actions are obtained
similarly. In particular,
\[
L_3=L_1^{-1},
\qquad
L_2=L_1^2L_0.
\]

Setting $a=L_1$ and $b=L_0$,
we have
\[
a^4=b^2=1,
\qquad
bab=a^{-1}.
\]
Thus the permutation group is
\[
F=\langle a,b\rangle\cong D_8,
\]
where $D_8$ denotes the dihedral
group of order $8$.

Since $M$ is a monoid of I-type,
its group of fractions fits into
the exact sequence
\[
1\longrightarrow\mathbb Z^4
\longrightarrow G(M)
\longrightarrow D_8
\longrightarrow1.
\]
In particular, $G(M)$ is a solvable
Bieberbach group with Hirsch length $4$.
\end{example}

\subsection{Comparison of the permutation groups}

Table~\ref{tab:permutation-groups} summarizes
the permutation groups associated with
the $23$ algebras listed in
Appendix~\ref{app:binomial-algebras}.

\begin{table}[htbp]
\centering
\renewcommand{\arraystretch}{1.2}

\begin{tabular}{|l|l|c|}
\hline
Permutation group $F$ & Algebras & Count\\
\hline

$\{1\}$ &
1 & 1\\

$C_2$ &
2, 4, 6, 7, 8, 10, 16 & 7\\

$C_3$ &
3, 18, 19, 20 & 4\\

$C_4$ &
17, 23 & 2\\

$C_2\times C_2$ &
5, 9, 11, 13, 14, 15, 22 & 7\\

$D_8$ &
12, 21 & 2\\

\hline
\end{tabular}

\vspace{0.3cm}
\caption{Permutation groups associated with
the four-dimensional binomial algebras.}
\label{tab:permutation-groups}
\end{table}

The permutation groups in
Table~\ref{tab:permutation-groups}
are all finite and solvable.
For the corresponding monoids of I-type,
their groups of fractions are
torsion-free, virtually abelian
Bieberbach groups of Hirsch length $4$.

\section{Computing the cocycles}\label{sec:computationscocycles}

We now study the second cohomology groups
of the groups of fractions considered
in the previous sections.

Throughout this section, we regard
$k^\times$ as a trivial $G(M)$-module.
We compute $H^2(G(M),k^\times)$
using group homology and the
universal coefficient theorem.

Recall that for a group $G$, the
universal coefficient theorem
\cite[Theorem~3.6.5]{weibel1994}
gives a short exact sequence
\[
\begin{aligned}
0\longrightarrow\;&
\operatorname{Ext}^1_{\mathbb Z}
(H_1(G,\mathbb Z),k^\times)\\
\longrightarrow\;&
H^2(G,k^\times)\\
\longrightarrow\;&
\operatorname{Hom}
(H_2(G,\mathbb Z),k^\times)
\longrightarrow 0.
\end{aligned}
\]
In particular, when $H_1(G,\mathbb Z)$
is free abelian, the Ext term vanishes.

\subsection{Dimension \texorpdfstring{$2$}{2}}

We begin with the two groups
arising in dimension $2$.

\begin{example}
For the polynomial ring in two variables,
the group of fractions is
\[
G(M)\cong\mathbb Z^2.
\]
Since
\[
H_2(\mathbb Z^n,\mathbb Z)
\cong\bigwedge^2\mathbb Z^n,
\]
and
\[
\bigwedge^2\mathbb Z^2\cong\mathbb Z,
\]
the universal coefficient theorem gives
\[
\begin{aligned}
H^2(G(M),k^\times)
&\cong
\operatorname{Hom}(\mathbb Z,k^\times)\\
&\cong k^\times.
\end{aligned}
\]
\end{example}

\begin{example}
For the coefficient-free relation
$x^2=y^2$, we identified the group
of fractions with the Klein bottle group
\[
\Gamma=
\langle a,b\mid aba^{-1}=b^{-1}\rangle.
\]

Since the Klein bottle is aspherical,
its cohomology agrees with the group
cohomology of $\Gamma$. We have
\[
H_1(\Gamma,\mathbb Z)
\cong\mathbb Z\oplus\mathbb Z/2\mathbb Z,
\qquad
H_2(\Gamma,\mathbb Z)=0.
\]
By the universal coefficient theorem,
\[
\begin{aligned}
H^2(\Gamma,k^\times)
&\cong
\operatorname{Ext}^1_{\mathbb Z}
(H_1(\Gamma,\mathbb Z),k^\times)\\
&\cong k^\times/(k^\times)^2.
\end{aligned}
\]

In particular, if $k$ is algebraically
closed, every nonzero element is a square,
so
\[
H^2(\Gamma,k^\times)=0.
\]
\end{example}

\subsection{Dimension \texorpdfstring{$3$}{3}}

We next consider the polynomial
ring in three variables.

\begin{example}
For the polynomial ring in three
variables, we have
\[
G(M)\cong\mathbb Z^3.
\]

Since
\[
H_2(\mathbb Z^3,\mathbb Z)
\cong\bigwedge^2\mathbb Z^3
\cong\mathbb Z^3,
\]
the universal coefficient theorem gives
\[
\begin{aligned}
H^2(G(M),k^\times)
&\cong
\operatorname{Hom}(\mathbb Z^3,k^\times)\\
&\cong(k^\times)^3.
\end{aligned}
\]
\end{example}

\subsection{Dimension \texorpdfstring{$4$}{4}}

We now turn to the four-dimensional examples
from the previous section. For each of the
$23$ algebras, we consider the corresponding
group of fractions $G(M)$ and compute
its first and second integral homology groups,
together with $H^2(G(M),k^\times)$.

We begin with a few examples to illustrate
the methods used in these computations.
The results for all $23$ algebras are
summarized in Table~\ref{tab:dimension4-cohomology}.

\begin{example}[Algebra 1]
For Algebra 1, the group of fractions is
\[
G(M)\cong\mathbb Z^4.
\]
For a free abelian group, we have
\[
H_q(\mathbb Z^n,\mathbb Z)
\cong\bigwedge^q\mathbb Z^n.
\]
In particular,
\[
H_1(G(M),\mathbb Z)\cong\mathbb Z^4,
\qquad
H_2(G(M),\mathbb Z)
\cong\bigwedge^2\mathbb Z^4
\cong\mathbb Z^6.
\]

Since $H_1(G(M),\mathbb Z)$ is free abelian,
the Ext term in the universal coefficient
theorem vanishes. Therefore,
\[
\begin{aligned}
H^2(G(M),k^\times)
&\cong
\operatorname{Hom}
(H_2(G(M),\mathbb Z),k^\times)\\
&\cong
\operatorname{Hom}(\mathbb Z^6,k^\times)\\
&\cong(k^\times)^6.
\end{aligned}
\]
\end{example}

\begin{example}[Algebra 2]
For Algebra 2, the group of fractions
has presentation
\[
G(M)=
\left\langle x_0,x_1,x_2,x_3
\ \middle|\
\begin{aligned}
x_0x_1&=x_1x_0,\\
x_0x_2&=x_2x_0,\\
x_0x_3&=x_3x_0,\\
x_1x_2&=x_2x_1,\\
x_1x_3&=x_3x_2,\\
x_2x_3&=x_3x_1
\end{aligned}
\right\rangle.
\]

The first three relations show that
$x_0$ is central. The remaining
relations imply
\[
x_1x_2=x_2x_1,\qquad
x_3^{-1}x_1x_3=x_2,\qquad
x_3^{-1}x_2x_3=x_1.
\]

Consequently,
\[
G(M)\cong
\left(\mathbb Z^2
\rtimes_\psi\mathbb Z\right)
\times\mathbb Z,
\qquad
\psi=
\begin{pmatrix}
0&1\\
1&0
\end{pmatrix}.
\]

Let
\[
K=\mathbb Z^2\rtimes_\psi\mathbb Z.
\]
Since
\[
\operatorname{coker}(\psi-I)
\cong\mathbb Z,
\]
we have
\[
H_1(K,\mathbb Z)\cong\mathbb Z^2.
\]
Thus
\[
H_1(G(M),\mathbb Z)\cong\mathbb Z^3.
\]

To compute the second homology of $K$,
we use the long exact homology sequence
associated with the mapping torus
of the torus automorphism $\psi$
\cite[Example~2.48]{hatcher_algebraic_2001}:
\[
\begin{aligned}
0\longrightarrow\;&
\operatorname{coker}
(\bigwedge^2\psi-I)\\
\longrightarrow\;&
H_2(K,\mathbb Z)\\
\longrightarrow\;&
\ker(\psi-I)
\longrightarrow0.
\end{aligned}
\]

Since $\det(\psi)=-1$, the induced
action on
$\bigwedge^2\mathbb Z^2\cong\mathbb Z$
is multiplication by $-1$. Hence
\[
\operatorname{coker}
(\bigwedge^2\psi-I)
\cong\mathbb Z/2\mathbb Z.
\]
Moreover,
\[
\ker(\psi-I)\cong\mathbb Z.
\]
The short exact sequence splits,
so
\[
H_2(K,\mathbb Z)
\cong\mathbb Z\oplus\mathbb Z/2\mathbb Z.
\]

By the K\"unneth formula
\cite[Theorem~3B.6]{hatcher_algebraic_2001},
\[
\begin{aligned}
H_2(G(M),\mathbb Z)
&\cong H_2(K,\mathbb Z)
\oplus H_1(K,\mathbb Z)\\
&\cong\mathbb Z^3
\oplus\mathbb Z/2\mathbb Z.
\end{aligned}
\]

Since $H_1(G(M),\mathbb Z)$
is free abelian, the Ext term
in the universal coefficient
theorem vanishes. Therefore,
\[
\begin{aligned}
H^2(G(M),k^\times)
&\cong
\operatorname{Hom}
(H_2(G(M),\mathbb Z),k^\times)\\
&\cong
(k^\times)^3\times\mu_2(k),
\end{aligned}
\]
where
\[
\mu_2(k)=
\{\lambda\in k^\times:\lambda^2=1\}.
\]

In particular, if
$\operatorname{char}(k)\ne2$,
then
\[
H^2(G(M),k^\times)
\cong(k^\times)^3\times C_2.
\]
\end{example}

\begin{example}[Algebra 3]
For Algebra 3, the group of fractions
has presentation
\[
G(M)=
\left\langle x_0,x_1,x_2,x_3
\ \middle|\
\begin{aligned}
x_0x_1&=x_1x_0,\\
x_0x_2&=x_2x_0,\\
x_0x_3&=x_3x_2,\\
x_1x_2&=x_2x_1,\\
x_1x_3&=x_3x_0,\\
x_2x_3&=x_3x_1
\end{aligned}
\right\rangle.
\]

The generators $x_0,x_1,x_2$ commute,
and the remaining relations give
\[
x_3^{-1}x_0x_3=x_2,\qquad
x_3^{-1}x_1x_3=x_0,\qquad
x_3^{-1}x_2x_3=x_1.
\]
Consequently,
\[
G(M)\cong
\mathbb Z^3\rtimes_\psi\mathbb Z,
\qquad
\psi=
\begin{pmatrix}
0&1&0\\
0&0&1\\
1&0&0
\end{pmatrix}.
\]

Since $\psi^3=I$, the subgroup
\[
\langle x_0,x_1,x_2,x_3^3\rangle
\cong\mathbb Z^4
\]
has index $3$ in $G(M)$.

The first homology group is
\[
H_1(G(M),\mathbb Z)
\cong
\mathbb Z\oplus
\operatorname{coker}(\psi-I).
\]
Since
\[
\operatorname{coker}(\psi-I)
\cong\mathbb Z,
\]
we obtain
\[
H_1(G(M),\mathbb Z)
\cong\mathbb Z^2.
\]

To compute the second homology,
we use the long exact homology sequence
associated with the mapping torus
of the torus automorphism $\psi$
\cite[Example~2.48]{hatcher_algebraic_2001}:
\[
\begin{aligned}
0\longrightarrow\;&
\operatorname{coker}
(\bigwedge^2\psi-I)\\
\longrightarrow\;&
H_2(G(M),\mathbb Z)\\
\longrightarrow\;&
\ker(\psi-I)
\longrightarrow0.
\end{aligned}
\]

With respect to the basis
\[
e_0\wedge e_1,\quad
e_0\wedge e_2,\quad
e_1\wedge e_2,
\]
the induced action is
\[
\bigwedge^2\psi=
\begin{pmatrix}
0&0&1\\
-1&0&0\\
0&-1&0
\end{pmatrix}.
\]
Both $\psi-I$ and
$\bigwedge^2\psi-I$ have Smith
normal form
\[
\operatorname{diag}(1,1,0).
\]
Therefore,
\[
\operatorname{coker}
(\bigwedge^2\psi-I)
\cong\mathbb Z,
\qquad
\ker(\psi-I)\cong\mathbb Z.
\]

The sequence splits, giving
\[
H_2(G(M),\mathbb Z)
\cong\mathbb Z^2.
\]

Since $H_1(G(M),\mathbb Z)$
is free abelian, the Ext term
in the universal coefficient
theorem vanishes. Hence
\[
\begin{aligned}
H^2(G(M),k^\times)
&\cong
\operatorname{Hom}
(H_2(G(M),\mathbb Z),k^\times)\\
&\cong
\operatorname{Hom}(\mathbb Z^2,k^\times)\\
&\cong(k^\times)^2.
\end{aligned}
\]
\end{example}

We carry out similar computations for all
$23$ four-dimensional binomial algebras
listed in Appendix~\ref{app:binomial-algebras}.
The resulting homology and cohomology
groups are summarized in
Table~\ref{tab:dimension4-cohomology}.

\begin{table}[H]
\centering

\scriptsize
\renewcommand{\arraystretch}{1.19}
\setlength{\tabcolsep}{5pt}

\begin{tabular}{|c|c|l|l|l|}
\hline
Algebra & $F$ &
$H_1(G(M),\mathbb Z)$ &
$H_2(G(M),\mathbb Z)$ &
$H^2(G(M),k^\times)$\\
\hline\hline
1  & $1$
   & $\mathbb Z^4$
   & $\mathbb Z^6$
   & $(k^\times)^6$\\ \hline

2  & $C_2$
   & $\mathbb Z^3$
   & $\mathbb Z^3\oplus C_2$
   & $(k^\times)^3\times C_2$\\ \hline

3  & $C_3$
   & $\mathbb Z^2$
   & $\mathbb Z^2$
   & $(k^\times)^2$\\ \hline

4  & $C_2$
   & $\mathbb Z^3$
   & $\mathbb Z^3\oplus C_2$
   & $(k^\times)^3\times C_2$\\ \hline

5  & $C_2^2$
   & $\mathbb Z^2\oplus C_2$
   & $\mathbb Z\oplus C_2\oplus C_4$
   & $k^\times\times C_2\times C_4$\\ \hline

6  & $C_2$
   & $\mathbb Z^3\oplus C_2$
   & $\mathbb Z^3\oplus C_2^2$
   & $(k^\times)^3\times C_2^2$\\ \hline

7  & $C_2$
   & $\mathbb Z^3$
   & $\mathbb Z^3\oplus C_2$
   & $(k^\times)^3\times C_2$\\ \hline

8  & $C_2$
   & $\mathbb Z^3$
   & $\mathbb Z^3\oplus C_2$
   & $(k^\times)^3\times C_2$\\ \hline

9  & $C_2^2$
   & $\mathbb Z^2\oplus C_2$
   & $\mathbb Z\oplus C_2^2$
   & $k^\times\times C_2^2$\\ \hline

10 & $C_2$
   & $\mathbb Z^2\oplus C_2$
   & $\mathbb Z^2\oplus C_2$
   & $(k^\times)^2\times C_2$\\ \hline

11 & $C_2^2$
   & $\mathbb Z^2\oplus C_2$
   & $\mathbb Z\oplus C_2\oplus C_4$
   & $k^\times\times C_2\times C_4$\\ \hline

12 & $D_8$
   & $\mathbb Z\oplus C_2$
   & $C_2$
   & $C_2$\\ \hline

13 & $C_2^2$
   & $\mathbb Z^2\oplus C_2^2$
   & $\mathbb Z\oplus C_2^3$
   & $k^\times\times C_2^3$\\ \hline

14 & $C_2^2$
   & $\mathbb Z^2\oplus C_2$
   & $\mathbb Z\oplus C_2^2$
   & $k^\times\times C_2^2$\\ \hline

15 & $C_2^2$
   & $\mathbb Z^2\oplus C_2$
   & $\mathbb Z\oplus C_2\oplus C_4$
   & $k^\times\times C_2\times C_4$\\ \hline

16 & $C_2$
   & $\mathbb Z^2\oplus C_2$
   & $\mathbb Z^2\oplus C_2$
   & $(k^\times)^2\times C_2$\\ \hline

17 & $C_4$
   & $\mathbb Z\oplus C_2^2$
   & $\mathbb Z\oplus C_2$
   & $k^\times\times C_2$\\ \hline

18 & $C_3$
   & $\mathbb Z^2\oplus C_3$
   & $\mathbb Z^2\oplus C_3$
   & $(k^\times)^2\times C_3$\\ \hline

19 & $C_3$
   & $\mathbb Z^2$
   & $\mathbb Z^2$
   & $(k^\times)^2$\\ \hline

20 & $C_3$
   & $\mathbb Z^2$
   & $\mathbb Z^2$
   & $(k^\times)^2$\\ \hline

21 & $D_8$
   & $\mathbb Z\oplus C_4$
   & $C_4$
   & $C_4$\\ \hline

22 & $C_2^2$
   & $\mathbb Z\oplus C_4$
   & $C_4$
   & $C_4$\\ \hline

23 & $C_4$
   & $\mathbb Z\oplus C_4$
   & $\mathbb Z$
   & $k^\times$\\ \hline
\end{tabular}
\vspace{0.3cm}
\caption{Homology and Cohomology
of the groups of fractions in dimension four.}
\label{tab:dimension4-cohomology}
\end{table}

\section{Explicit cocycle twists of some algebras}
\label{sec:explicit-cocycle-twists}

In Section~\ref{sec:computationscocycles}, we computed
$H^2(G(M),k^\times)$ for the $23$ binomial algebras.
We now use these computations to construct families
of algebras with nonzero coefficients in their
defining relations.

As explained in
Section~\ref{binomasregulaandtwistedmonoid},
a group $2$-cocycle determines the coefficients
of the binomial relations in a twisted monoid
algebra. We illustrate this construction for
Algebras 1, 2 and 3.

\begin{example}[Algebra 1]
For Algebra 1, we have
\[
G(M)\cong\mathbb Z^4,
\qquad
H^2(G(M),k^\times)\cong(k^\times)^6.
\]

Let $e_0,e_1,e_2,e_3$ be the standard basis
of $\mathbb Z^4$. Given a group $2$-cocycle
$\sigma$, set
\[
q_{ij}
=
\frac{\sigma(e_i,e_j)}
{\sigma(e_j,e_i)},
\qquad 0\leq i<j\leq3.
\]
These six scalars determine the cohomology
class of $\sigma$, and every choice of
nonzero scalars $q_{ij}$ can be obtained
from a group $2$-cocycle.

The resulting twisted algebra has relations
\[
\begin{aligned}
x_0x_1&=q_{01}x_1x_0,
&
x_0x_2&=q_{02}x_2x_0,\\
x_0x_3&=q_{03}x_3x_0,
&
x_1x_2&=q_{12}x_2x_1,\\
x_1x_3&=q_{13}x_3x_1,
&
x_2x_3&=q_{23}x_3x_2.
\end{aligned}
\]
Here the six parameters are arbitrary elements
of $k^\times$. This gives the family of
multiparameter quantum affine spaces.
Setting all parameters equal to $1$
recovers Algebra 1.
\end{example}

\begin{example}[Algebra 2]

For Algebra 2, we obtained
\[
G(M)\cong
(\mathbb Z^2\rtimes_\psi\mathbb Z)
\times\mathbb Z,
\]
where $\psi$ exchanges the two standard
generators of $\mathbb Z^2$. We also computed
\[
H^2(G(M),k^\times)
\cong(k^\times)^3\times C_2.
\]

To construct cocycles explicitly, consider
three invertible generators $a,b,c$
satisfying
\[
ab=pba,\qquad
ac=pca,\qquad
bc=\varepsilon cb,
\]
where
\[
p\in k^\times,
\qquad
\varepsilon\in\{1,-1\}.
\]

We introduce another invertible generator
$d$ satisfying
\[
da=qad,\qquad
db=cd,\qquad
dc=rbd,
\]
where $q,r\in k^\times$.
These relations are compatible because
$\varepsilon^2=1$.
They give a twisted group algebra
of $G(M)$.

Restricting the cocycle to the monoid
gives the following family:
\[
\begin{aligned}
x_0x_1&=p\,x_1x_0,
&
x_0x_2&=p\,x_2x_0,\\
x_0x_3&=q\,x_3x_0,
&
x_1x_2&=\varepsilon\,x_2x_1,\\
x_1x_3&=r^{-1}x_3x_2,
&
x_2x_3&=x_3x_1.
\end{aligned}
\]
Here
\[
p,q,r\in k^\times,
\qquad
\varepsilon\in\{1,-1\}.
\]

The three nonzero parameters correspond
to the factors $(k^\times)^3$, while
$\varepsilon$ comes from the order-two
factor of the cohomology group.
In particular, the coefficients of
the first two relations must agree,
and the coefficient of the fourth
relation must satisfy
$\varepsilon^2=1$.

Setting $p=q=r=\varepsilon=1$
recovers Algebra 2.
\end{example}

\begin{example}[Algebra 3]

For Algebra 3, we obtained
\[
G(M)\cong
\mathbb Z^3\rtimes_\psi\mathbb Z,
\]
where $\psi$ cyclically permutes the
three standard generators of $\mathbb Z^3$.
We also computed
\[
H^2(G(M),k^\times)
\cong(k^\times)^2.
\]

To construct cocycles explicitly, consider
three invertible generators $a,b,c$
satisfying
\[
ab=tba,\qquad
ac=t^{-1}ca,\qquad
bc=tcb,
\]
where $t\in k^\times$.

We introduce another invertible generator
$d$ satisfying
\[
ad=s^{-1}dc,\qquad
bd=da,\qquad
cd=db,
\]
where $s\in k^\times$.
These relations are compatible with
the relations among $a,b,c$ and
describe a cocycle twist of the
group algebra $k[G(M)]$.

Restricting the cocycle to the monoid
gives the following family:
$$ \begin{aligned}
x_0x_1&=t\,x_1x_0,
&
x_0x_2&=t^{-1}x_2x_0,\\
x_0x_3&=s^{-1}x_3x_2,
&
x_1x_2&=t\,x_2x_1,\\
x_1x_3&=x_3x_0,
&
x_2x_3&=x_3x_1.
\end{aligned}
$$
Here $t,s\in k^\times$ are arbitrary.

The cyclic action relates the coefficients
of the first, second and fourth relations.
In particular, the coefficients of the
first and fourth relations are equal,
while the coefficient of the second
relation is their inverse.

Thus we obtain a family with two
independent nonzero parameters,
in agreement with the computation
of $H^2(G(M),k^\times)$.

Setting $t=s=1$ recovers Algebra 3.
\end{example}

The constructions above give families of
four-dimensional Koszul Artin--Schelter regular
algebras through cocycle twists.
Similar computations can be carried out
for the remaining binomial algebras
to obtain further families.
In future work, we would like to compute
the Kodaira--Spencer maps and discrete invariants
of these families, following the methods in
\cite{bhatoy-colin-felix-ravali-components,bhatoy2026discreteinvariants}.
We would also like to investigate whether
these families belong to the components
of the moduli space studied in our previous work,
or whether they give rise to new families.
Another possible direction is to extend
these computations to dimension five,
using the classification of Yang--Baxter
solutions in~\cite[Table 2.1]{AkgunMerebVendramin2022}.

\clearpage

\appendix
\section{Defining relations of the binomial algebras}
\label{app:binomial-algebras}

The following tables list the defining relations of the four-dimensional
binomial Artin--Schelter regular algebras considered in this paper.\\

\begingroup
\footnotesize
\setlength{\tabcolsep}{1.5pt}
\renewcommand{\arraystretch}{0.9}
\setlength{\LTpre}{4pt}
\setlength{\LTpost}{4pt}
\setlength{\LTleft}{\fill}
\setlength{\LTright}{\fill}

\newcommand{\algentry}[2]{%
\begin{minipage}[t]{0.29\textwidth}
\centering
\textbf{Algebra #1}\\[2pt]
$\begin{gathered}
#2
\end{gathered}$
\end{minipage}%
}


\begin{longtable}{|c|c|c|}

\hline
\endfirsthead
\hline
\endhead
\endfoot
\endlastfoot


\algentry{1}{
x_0x_1-x_1x_0=0\\
x_0x_2-x_2x_0=0\\
x_0x_3-x_3x_0=0\\
x_1x_2-x_2x_1=0\\
x_1x_3-x_3x_1=0\\
x_2x_3-x_3x_2=0
}
&
\algentry{2}{
x_0x_1-x_1x_0=0\\
x_0x_2-x_2x_0=0\\
x_0x_3-x_3x_0=0\\
x_1x_2-x_2x_1=0\\
x_1x_3-x_3x_2=0\\
x_2x_3-x_3x_1=0
}
&
\algentry{3}{
x_0x_1-x_1x_0=0\\
x_0x_2-x_2x_0=0\\
x_0x_3-x_3x_2=0\\
x_1x_2-x_2x_1=0\\
x_1x_3-x_3x_0=0\\
x_2x_3-x_3x_1=0
}\\[2pt]
\hline


\algentry{4}{
x_0x_1-x_1x_0=0\\
x_0x_2-x_2x_1=0\\
x_0x_3-x_3x_1=0\\
x_1x_2-x_2x_0=0\\
x_1x_3-x_3x_0=0\\
x_2x_3-x_3x_2=0
}
&
\algentry{5}{
x_0x_1-x_1x_0=0\\
x_0x_2-x_3x_1=0\\
x_0x_3-x_2x_1=0\\
x_1x_2-x_3x_0=0\\
x_1x_3-x_2x_0=0\\
x_2x_3-x_3x_2=0
}
&
\algentry{6}{
x_0^2-x_1^2=0\\
x_0x_2-x_2x_0=0\\
x_0x_3-x_3x_0=0\\
x_1x_2-x_2x_1=0\\
x_1x_3-x_3x_1=0\\
x_2x_3-x_3x_2=0
}\\[2pt]
\hline


\algentry{7}{
x_0^2-x_1^2=0\\
x_0x_2-x_2x_0=0\\
x_0x_3-x_3x_1=0\\
x_1x_2-x_2x_1=0\\
x_1x_3-x_3x_0=0\\
x_2x_3-x_3x_2=0
}
&
\algentry{8}{
x_0^2-x_1^2=0\\
x_0x_2-x_2x_1=0\\
x_0x_3-x_3x_1=0\\
x_1x_2-x_2x_0=0\\
x_1x_3-x_3x_0=0\\
x_2x_3-x_3x_2=0
}
&
\algentry{9}{
x_0^2-x_1^2=0\\
x_0x_2-x_2x_0=0\\
x_0x_3-x_3x_0=0\\
x_1x_2-x_3x_1=0\\
x_1x_3-x_2x_1=0\\
x_2x_3-x_3x_2=0
}\\[2pt]
\hline


\algentry{10}{
x_0^2-x_1^2=0\\
x_0x_2-x_3x_0=0\\
x_0x_3-x_2x_0=0\\
x_1x_2-x_3x_1=0\\
x_1x_3-x_2x_1=0\\
x_2x_3-x_3x_2=0
}
&
\algentry{11}{
x_0^2-x_1^2=0\\
x_0x_2-x_3x_1=0\\
x_0x_3-x_2x_1=0\\
x_1x_2-x_3x_0=0\\
x_1x_3-x_2x_0=0\\
x_2x_3-x_3x_2=0
}
&
\algentry{12}{
x_0^2-x_2x_3=0\\
x_0x_2-x_3x_0=0\\
x_0x_3-x_1x_2=0\\
x_1^2-x_3x_2=0\\
x_1x_3-x_2x_1=0\\
x_2x_0-x_3x_1=0
}\\[2pt]
\hline


\algentry{13}{
x_0^2-x_1^2=0\\
x_0x_2-x_2x_0=0\\
x_0x_3-x_3x_0=0\\
x_1x_2-x_2x_1=0\\
x_1x_3-x_3x_1=0\\
x_2^2-x_3^2=0
}
&
\algentry{14}{
x_0^2-x_1^2=0\\
x_0x_2-x_2x_0=0\\
x_0x_3-x_3x_1=0\\
x_1x_2-x_2x_1=0\\
x_1x_3-x_3x_0=0\\
x_2^2-x_3^2=0
}
&
\algentry{15}{
x_0^2-x_1^2=0\\
x_0x_2-x_2x_1=0\\
x_0x_3-x_3x_1=0\\
x_1x_2-x_2x_0=0\\
x_1x_3-x_3x_0=0\\
x_2^2-x_3^2=0
}\\[2pt]
\hline
\noalign{\vskip 0.3cm}
\caption{Defining relations of binomial
Artin--Schelter regular algebras 1--15.}
\label{tab:binomial-algebras-1}\\

\end{longtable}
\endgroup

\begingroup
\footnotesize
\setlength{\tabcolsep}{1.5pt}
\renewcommand{\arraystretch}{0.9}
\setlength{\LTpre}{4pt}
\setlength{\LTpost}{4pt}
\setlength{\LTleft}{\fill}
\setlength{\LTright}{\fill}

\newcommand{\algentry}[2]{%
\begin{minipage}[t]{0.29\textwidth}
\centering
\textbf{Algebra #1}\\[2pt]
$\begin{gathered}
#2
\end{gathered}$
\end{minipage}%
}

\begin{longtable}{|c|c|c|}

\hline
\endfirsthead
\hline
\endhead
\endfoot
\endlastfoot

\algentry{16}{
x_0^2-x_1^2=0\\
x_0x_2-x_3x_1=0\\
x_0x_3-x_2x_1=0\\
x_1x_2-x_3x_0=0\\
x_1x_3-x_2x_0=0\\
x_2^2-x_3^2=0
}
&
\algentry{17}{
x_0^2-x_3^2=0\\
x_0x_2-x_1x_3=0\\
x_0x_3-x_2x_1=0\\
x_1^2-x_2^2=0\\
x_1x_2-x_3x_0=0\\
x_2x_0-x_3x_1=0
}
&
\algentry{18}{
x_0^2-x_1x_2=0\\
x_0x_1-x_2^2=0\\
x_0x_3-x_3x_0=0\\
x_1^2-x_2x_0=0\\
x_1x_3-x_3x_1=0\\
x_2x_3-x_3x_2=0
}\\[2pt]
\hline

\algentry{19}{
x_0^2-x_1x_2=0\\
x_0x_1-x_2^2=0\\
x_0x_3-x_3x_2=0\\
x_1^2-x_2x_0=0\\
x_1x_3-x_3x_0=0\\
x_2x_3-x_3x_1=0
}
&
\algentry{20}{
x_0^2-x_1x_2=0\\
x_0x_1-x_2^2=0\\
x_0x_3-x_3x_1=0\\
x_1^2-x_2x_0=0\\
x_1x_3-x_3x_2=0\\
x_2x_3-x_3x_0=0
}
&
\algentry{21}{
x_0^2-x_3x_1=0\\
x_0x_1-x_1x_2=0\\
x_0x_2-x_2x_0=0\\
x_1^2-x_3^2=0\\
x_1x_3-x_2^2=0\\
x_2x_3-x_3x_0=0
}\\[2pt]
\hline


\algentry{22}{
x_0^2-x_3x_1=0\\
x_0x_1-x_2x_3=0\\
x_0x_2-x_1^2=0\\
x_1x_2-x_3x_0=0\\
x_1x_3-x_2^2=0\\
x_2x_0-x_3^2=0
}
&
\multicolumn{2}{>{\centering\arraybackslash}p{\dimexpr0.62\textwidth+2\tabcolsep+\arrayrulewidth\relax}|}{%
\algentry{23}{
x_0^2-x_1x_3=0\\
x_0x_1-x_2x_3=0\\
x_0x_2-x_3^2=0\\
x_1^2-x_2x_0=0\\
x_1x_2-x_3x_0=0\\
x_2^2-x_3x_1=0
}
}\\[2pt]
\hline
\noalign{\vskip 0.3cm}
\caption{Defining relations of binomial
Artin--Schelter regular algebras 16--23.}
\label{tab:binomial-algebras-2}\\

\end{longtable}

\endgroup

\pagebreak 

\bibliographystyle{plain}
\bibliography{bibliography}

\end{document}